\documentclass[12pt]{article}
\usepackage{amssymb,amsmath}
\usepackage{hyperref}
\usepackage{graphicx}
\usepackage{color}

\begin{document}

\title{\LARGE\bf A rational approximation of the arctangent function and a new approach in computing pi}

\author{
\normalsize\bf S. M. Abrarov\footnote{\scriptsize{Dept. Earth and Space Science and Engineering, York University, Toronto, Canada, M3J 1P3.}}\, and B. M. Quine$^{*}$\footnote{\scriptsize{Dept. Physics and Astronomy, York University, Toronto, Canada, M3J 1P3.}}}

\date{March 24, 2016}
\maketitle

\begin{abstract}
We have shown recently that integration of the error function ${\rm{erf}}\left( x \right)$ represented in form of a sum of the Gaussian functions provides an asymptotic expansion series for the constant pi. In this work we derive a rational approximation of the arctangent function $\arctan \left( x \right)$ that can be readily generalized it to its counterpart $ - {\rm{sgn}}\left( x \right)\pi /2 + \arctan \left( x \right)$, where ${\rm{sgn}}\left( x \right)$ is the signum function. The application of the expansion series for these two functions leads to a new asymptotic formula for $\pi$.
\vspace{0.25cm}
\vspace{0.25cm}
\end{abstract}

\section{Derivation}

Consider the following integral \cite{Fayed2014}
\begin{equation}\label{eq_1}
\int\limits_0^\infty  {{e^{ - {y ^2}{t^2}}}{\rm{erf}}\left( {x t} \right)dt}  = \frac{1}{{y \sqrt \pi  }}\arctan \left( {\frac{x }{y }} \right),
\end{equation}
\\
\noindent where we imply that all variables $t$, $x$ and $y$ are real. Assuming that $y  = 1$ the integral \eqref{eq_1} can be rewritten as
\begin{equation}\label{eq_2}
\arctan \left( x \right) = \sqrt \pi  \int\limits_0^\infty  {{e^{ - {t^2}}}{\rm{erf}}\left( {xt} \right)dt}.
\end{equation}

The error function can be represented in form of a sum of the Gaussian functions (see Appendix A)
\begin{equation}\label{eq_3}
{\rm{erf}}\left( x \right) = \frac{{2x}}{{\sqrt \pi  }} \times \mathop {\lim }\limits_{L \to \infty } \frac{1}{L}\sum\limits_{\ell  = 1}^L {{e^{ - \frac{{{{\left( {\ell  - 1/2} \right)}^2}{x^2}}}{{{L^2}}}}}}.
\end{equation}
Consequently, substituting this limit into the equation \eqref{eq_2} leads to
\[
\arctan \left( x \right) = \sqrt \pi   \times \mathop {\lim }\limits_{L \to \infty } \int\limits_0^\infty  {{e^{ - {t^2}}}\underbrace {\frac{{2xt}}{{\sqrt \pi  L}}\sum\limits_{\ell  = 1}^L {{e^{ - \frac{{{{\left( {\ell  - 1/2} \right)}^2}{x^2}{t^2}}}{{{L^2}}}}}} }_{{\rm{erf}}\left( {xt} \right)}dt}.
\]
Each integral term in this equation is analytically integrable. Consequently, we obtain a new equation for the arctangent function
\begin{equation}\label{eq_4}
\arctan \left( x \right) = 4 \times \mathop {\lim }\limits_{L \to \infty } \sum\limits_{\ell  = 1}^L {\frac{{Lx}}{{{{\left( {2\ell  - 1} \right)}^2}{x^2} + 4{L^2}}}}.
\end{equation}

Since
$$
\pi  = 4\arctan \left( 1 \right)
$$
it follows that
\begin{equation}\label{eq_5}
\pi  = 16 \times \mathop {\lim }\limits_{L \to \infty } \sum\limits_{\ell  = 1}^L {\frac{L}{{{{\left( {2\ell  - 1} \right)}^2} + 4{L^2}}}}.
\end{equation}
It should be noted that the limit \eqref{eq_5} has been reported already in our recent work \cite{Abrarov2016}. 

\begin{figure}[ht]
\begin{center}
\includegraphics[width=22pc]{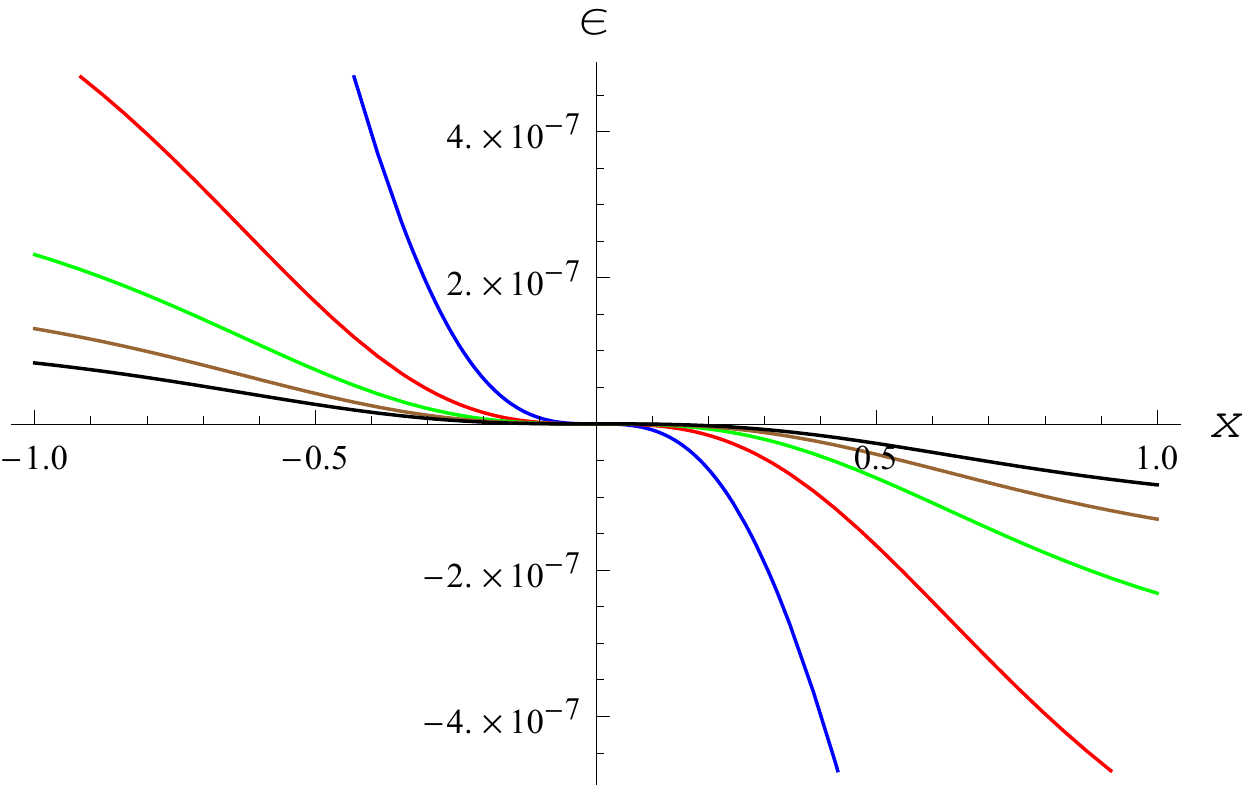}\hspace{2pc}%
\begin{minipage}[b]{28pc}
\vspace{0.5cm}
{\sffamily {\bf{Fig. 1.}} The difference $\varepsilon $ over the range $ - 1 \le x \le 1$ at $L = 100$ (blue), $L = 200$ (red), $L = 300$ (green), $L = 400$ (brown) and $L = 500$ (black).}
\end{minipage}
\end{center}
\end{figure}

\newpage
Truncation of the limit \eqref{eq_4} yields a rational approximation of the arctangent function
\begin{equation}\label{eq_6}
\arctan \left( x \right) \approx 4L\sum\limits_{\ell  = 1}^L {\frac{x}{{{{\left( {2\ell  - 1} \right)}^2}{x^2} + 4{L^2}}}}.
\end{equation}
Figure 1 shows the difference between the original arctangent function $\arctan \left( x \right)$ and its rational approximation \eqref{eq_6}
$$
\varepsilon  = \arctan \left( x \right) - 4L\sum\limits_{\ell  = 1}^L {\frac{x}{{{{\left( {2\ell  - 1} \right)}^2}{x^2} + 4{L^2}}}}
$$
over the range $ - 1 \le x \le 1$ at $L = 100$, $L = 200$, $L = 300$, $L = 400$ and $L = 500$ shown by blue, red, green, brown and black curves, respectively. As we can see from this figure, the difference $\varepsilon $ is dependent upon $x$. In particular, it increases with increasing argument by absolute value $\left| x \right|$. Thus, we can conclude that the rational approximation \eqref{eq_6} of the arctangent function is more accurate when its argument is smaller. Consequently, in order to obtain a higher accuracy we have to look for an equation in the form
$$
\pi  = \sum\limits_{n = 1}^N {{a_n}\arctan \left( {{b_n}} \right)},		\qquad \left|{b_n}\right| <  < 1,
$$
where ${a_n}$ and ${b_n}$ are the coefficients, with arguments of the arctangent function as small as possible by absolute value $\left|{b_n}\right|$. For example, applying the equation \eqref{eq_6} we may expect that at some fixed $L$ the approximation
$$
\pi = \, 4\arctan \left( {x = 1} \right) \approx 16L\sum\limits_{\ell  = 1}^L {\frac{1}{{{{\left( {2\ell  - 1} \right)}^2} + 4{L^2}}}}
$$
is less accurate than the approximation based on the Machin\text{'}s formula \cite{Borwein1989, Borwein2015}
\small
\[
\begin{aligned}
\pi  &= \, 4\left[ {4\arctan \left( {\frac{1}{5}} \right) - \arctan \left( {\frac{1}{{239}}} \right)} \right]\\
 &\approx 16L\sum\limits_{\ell  = 1}^L {\left( {\frac{{4\left( {1/5} \right)}}{{{{\left( {2\ell  - 1} \right)}^2}{{\left( {1/5} \right)}^2} + 4{L^2}}} - \frac{{1/239}}{{{{\left( {2\ell  - 1} \right)}^2}{{\left( {1/239} \right)}^2} + 4{L^2}}}} \right)}.
\end{aligned}
\]
\normalsize
Furthermore, with same equation \eqref{eq_6} for $\arctan \left( x \right)$ we can improve accuracy by using another formula for pi \cite{Borwein2015}
\small
\[
\begin{aligned}
\pi  = & \, 4\left[ {12\arctan \left( {\frac{1}{{18}}} \right) + 8\arctan \left( {\frac{1}{{57}}} \right) - 5\arctan \left( {\frac{1}{{239}}} \right)} \right]\\
 \approx & \, 16L\sum\limits_{\ell  = 1}^L \left( \frac{{12\left( {1/18} \right)}}{{{{\left( {2\ell  - 1} \right)}^2}{{\left( {1/18} \right)}^2} + 4{L^2}}} + \frac{{8\left( {1/57} \right)}}{{{{\left( {2\ell  - 1} \right)}^2}{{\left( {1/57} \right)}^2} + 4{L^2}}} \right. \\
& \left. - \frac{{5\left( {1/239} \right)}}{{{{\left( {2\ell  - 1} \right)}^2}{{\left( {1/239} \right)}^2} + 4{L^2}}} \right)
\end{aligned}
\]
\normalsize
due to smaller arguments $b_{n}$ of the arctangent function.

\section{Application}

\subsection{Counterpart function}

Once the rational approximation \eqref{eq_6} for the arctangent function is found, from the identity
$$
\arctan \left( {\frac{1}{x}} \right) + \arctan \left( x \right) = {\frac{\pi}{2}} {\rm{sgn}} \left( x \right),
$$
where
\[
{\rm{sgn}}\left( x \right) = \left\{ \begin{aligned}
1, \qquad & x > 0\\
0, \qquad & x = 0\\
 - 1, \qquad & x < 0
\end{aligned} \right.
\]
is the signum function \cite{Weisstein2003}, it follows that
\begin{equation}\label{eq_7}
- 4L\sum\limits_{\ell  = 1}^L {\frac{x}{{{{\left( {2\ell  - 1} \right)}^2} + 4{L^2}{x^2}}}} \approx - \frac{\pi }{2}{\rm{sgn}}\left( x \right) + \arctan \left( x \right).
\end{equation}

Figure 2 shows the expansion series \eqref{eq_7} computed at $L = 100$ (blue curve). The arctangent function is also shown for comparison (red curve). As we can see from this figure, on the left-half plane the expansion series \eqref{eq_7} is greater than the original arctangent function by $\pi /2$, while on the right-half plane it is smaller than the original arctangent function by $\pi /2$. 

\begin{figure}[ht]
\begin{center}
\includegraphics[width=22pc]{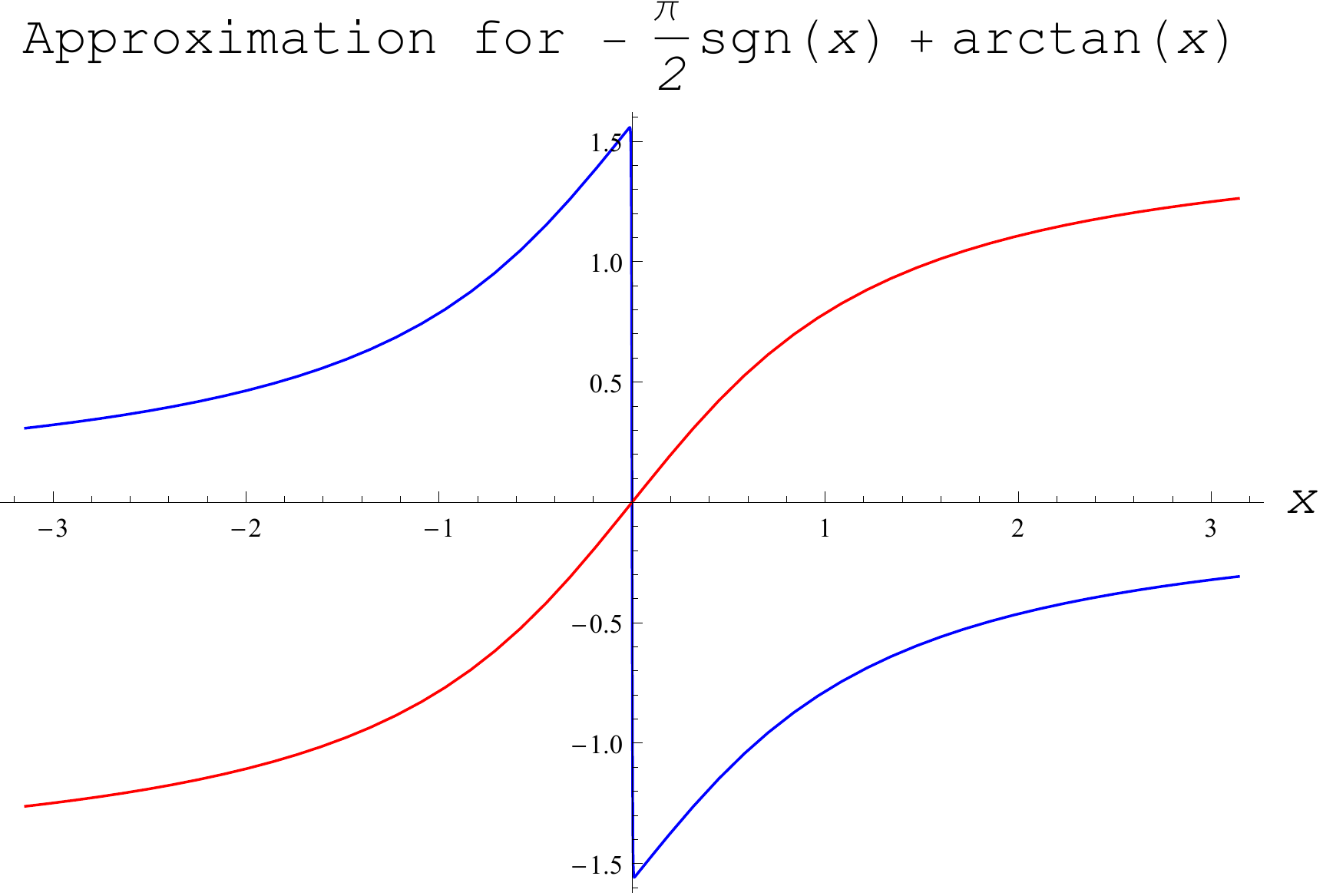}\hspace{2pc}%
\begin{minipage}[b]{28pc}
\vspace{0.75cm}
{\sffamily {\bf{Fig. 2.}} The expansion series \eqref{eq_7} computed at $L = 100$ (blue curve) resembling the function $ - {\rm{sgn}}\left( x \right)\pi /2 + \arctan \left( x \right)$. The original arctangent function (red curve) is also shown for comparison.}
\end{minipage}
\end{center}
\end{figure}

The approximation \eqref{eq_7} can be replaced with exact relation by tending the integer $L$ to infinity and taking the limit as
\begin{equation}\label{eq_8}
- 4 \times \mathop {\lim }\limits_{L \to \infty } \sum\limits_{\ell  = 1}^L {\frac{{Lx}}{{{{\left( {2\ell  - 1} \right)}^2} + 4{L^2}{x^2}}}}  =  - \frac{\pi }{2}{\rm{sgn}}\left( x \right) + \arctan \left( x \right).
\end{equation}
Since this limit represents a simple generalization of the equation \eqref{eq_4}, the function $- {\rm{sgn}} \left( x \right)\pi/2 + \arctan \left( x \right)$ can be regarded as a counterpart to the arctangent function $\arctan \left( x \right)$.

\subsection{Asymptotic formula for pi}

Using the limits \eqref{eq_4} and \eqref{eq_8} for the arctangent function $\arctan \left( x \right)$ and its counterpart function $ - {\rm{sgn}}\left( x \right)\pi /2 + \arctan \left( x \right)$, we can readily obtain an asymptotic expansion series for pi. Let us rewrite the equation \eqref{eq_8} as follows
\begin{equation}\label{eq_9}
\arctan \left( x \right) = \frac{\pi }{2}{\rm{sgn}}\left( x \right) - 4 \times \mathop {\lim }\limits_{L \to \infty } \sum\limits_{\ell  = 1}^L {\frac{{Lx}}{{{{\left( {2\ell  - 1} \right)}^2} + 4{L^2}{x^2}}}}.
\end{equation}
The difference of the equations \eqref{eq_9} and \eqref{eq_4} yields
\footnotesize
\[
0 = \underbrace {\frac{\pi }{2}{\rm{sgn}}\left( x \right) - \left( {4 \times \mathop {\lim }\limits_{L \to \infty } \sum\limits_{\ell  = 1}^L {\frac{{Lx}}{{{{\left( {2\ell  - 1} \right)}^2} + 4{L^2}{x^2}}}} } \right)}_{\,{\rm{eq}}{\rm{.}}\,\,\left( 9 \right)} - \underbrace {\left( {4 \times \mathop {\lim }\limits_{L \to \infty } \sum\limits_{\ell  = 1}^L {\frac{{Lx}}{{{{\left( {2\ell  - 1} \right)}^2}{x^2} + 4{L^2}}}} } \right)}_{{\rm{eq}}{\rm{.}}\,\,\left( 4 \right)}
\]
\normalsize
or
$$
4 \times \mathop {\lim }\limits_{L \to \infty } \sum\limits_{\ell  = 1}^L \left[ {\frac{{Lx}}{{{{\left( {2\ell  - 1} \right)}^2} + 4{L^2}{x^2}}}}  + \frac{{Lx}}{{{{\left( {2\ell  - 1} \right)}^2}{x^2} + 4{L^2}}} \right] = \frac{\pi }{2}{\rm{sgn}}\left( x \right)
$$
or
\begin{equation}\label{eq_10}
\pi = 8 \times \mathop {\lim }\limits_{L \to \infty } \sum\limits_{\ell  = 1}^L {L\left| x \right|\left[ {\frac{1}{{{{\left( {2\ell  - 1} \right)}^2}{x^2} + 4{L^2}}} + \frac{1}{{{{\left( {2\ell  - 1} \right)}^2} + 4{L^2}{x^2}}}} \right]}
\end{equation}
since ${\rm{sgn}}\left( x \right) = x/ \left| x \right|$ \cite{Weisstein2003}. Obviously, the equation \eqref{eq_10} can be interpreted as
\[
 \pi = 2 \frac{\left| x \right|}{x}\left[ \arctan \left( {\frac{1}{x}} \right) + \arctan \left( x \right) \right].
\]

Remarkably, although the argument $x$ is still present in the limit \eqref{eq_10} this asymptotic expansion series remains, nevertheless, independent of $x$. This signifies that according to equation \eqref{eq_10} the constant $\pi $ can be computed at any real value of the argument $x \in \mathbb{R}$.

The limit \eqref{eq_10} can be truncated by an arbitrarily large value $L >  > 1$ as given by
\begin{equation}\label{eq_11}
\pi  \approx 8L\left| x \right|\sum\limits_{\ell  = 1}^L {\left[ {\frac{1}{{{{\left( {2\ell  - 1} \right)}^2}{x^2} + 4{L^2}}} + \frac{1}{{{{\left( {2\ell  - 1} \right)}^2} + 4{L^2}{x^2}}}} \right]}.
\end{equation}
We performed sample computations by using Wolfram Mathematica 9 in enhanced precision mode in order to visualize the number of coinciding digits with actual value of the constant pi
$$
3.1415926535897932384626433832795028841971693993751 \ldots \,\,.
$$
The sample computations show that accuracy of the approximation limit \eqref{eq_11} depends upon the two values $L$ and  $x$ (the dependence on the argument $x$ in the equation \eqref{eq_11} is due to truncation now). For example, at $L = {10^{12}}$ and $x = 1$, we get
$$
\underbrace {3.141592653589793238462643}_{25\,\,{\rm{coinciding}}\,\,{\rm{digits}}}46661283621753050273271 \ldots \,\,,
$$
while at  same $L = {10^{12}}$ but smaller $x = {10^{ - 9}}$, the result is
$$
\underbrace {3.14159265358979323846264338327950}_{33\,\,{\rm{coinciding}}\,\,{\rm{digits}}}305086383606604 \ldots \,\,.
$$
Comparing these approximated values with the actual value for the constant pi one can see that at $x = 1$ and $x = {10^{ - 9}}$ the quantity of coinciding digits are $25$ and $33$, respectively. It should be noted, however, that the argument $x$ cannot be taken arbitrarily small since its optimized value depends upon the chosen integer $L.$

\section{Conclusion}
We obtain an efficient rational approximation for the arctangent function $\arctan \left( x \right)$ that can be generalized to its counterpart function $- {\rm{sgn}}\left( x \right)\pi /2 + \arctan \left( x \right)$. The application of the expansion series of the arctangent function and its counterpart results in a new formula for $\pi$. The computational test we performed shows that the new asymptotic expansion series for pi may be rapid in convergence.

\section*{Acknowledgments}

This work is supported by National Research Council Canada, Thoth Technology Inc. and York University. The authors would like to thank Prof. H. Rosengren and Prof. L. Tournier for review and useful information.

\section*{Appendix A}

Consider an integral of the error function (see integral 12 on page 4 in \cite{Ng1969})
\[
{\rm{erf}}\left( x \right) = \frac{1}{\pi }\int\limits_0^\infty  {{e^{ - u}}\sin \left( {2x\sqrt u } \right)\frac{{du}}{u}}.
\]
This integral can be readily expressed through the sinc function
$$
\left\{ {\rm{sinc}} \left( x \ne 0 \right) = {\rm{sin}} \left( x \right) / x, \,\, {\rm{sinc}} \left( x = 0 \right) = 1 \right\}
$$
by making change of the variable $v = \sqrt u$ leading to
\[
\begin{aligned}
{\rm{erf}}\left( x \right) &= \frac{1}{\pi }\int\limits_0^\infty  {{e^{ - {v^2}}}\sin \left( {2xv} \right)\frac{{2vdv}}{{{v^2}}}}  = \frac{2}{\pi }\int\limits_0^\infty  {{e^{ - {v^2}}}\sin \left( {2xv} \right)\frac{{dv}}{v}} \\
 &= \frac{{4x}}{\pi }\int\limits_0^\infty  {{e^{ - {v^2}}}\sin \left( {2xv} \right)\frac{{dv}}{{2xv}}} 
\end{aligned}
\]
or
\[
{\rm{erf}}\left( x \right) = \frac{{4x}}{\pi }\int\limits_0^\infty  {{e^{ - {v^2}}}{\rm{sinc}}\left( {2xv} \right)dv}.
\]
The factor $2$ in the argument of the sinc function can be excluded by making change of the variable $t = 2v$ again. This provides
\[
{\rm{erf}}\left( x \right) = \frac{{4x}}{\pi }\int\limits_0^\infty  {{e^{ - {t^2}/4}}{\rm{sinc}}\left( {xt} \right)\frac{{dt}}{2}}
\]
or
\[\label{A.1}
\tag{A.1}
{\rm{erf}}\left( x \right) = \frac{{2x}}{\pi }\int\limits_0^\infty  {{e^{ - {t^2}/4}}{\rm{sinc}}\left( {xt} \right)dt}.
\]

As it has been shown in our recent publication, the sinc function can be expressed as given by \cite{Abrarov2015}
\[\label{A.2}
\tag{A.2}
{\rm{sinc}}\left( x \right) = \mathop {\lim }\limits_{L \to \infty } \frac{1}{L}\sum\limits_{\ell  = 1}^L {\cos \left( {\frac{{\ell  - 1/2}}{L}x} \right)}.
\]
From the following integral
\[\label{A.3}
\tag{A.3}
{\rm{sinc}}\left( x \right) = \int\limits_0^1 {\cos \left( {xu} \right)du}  = \frac{1}{x}\int\limits_0^x {\cos \left( t \right)dt}
\]
it is not difficult to see that the cosine expansion \eqref{A.2} of the sinc function is just a result of integration of equation \eqref{A.3} performed by using the midpoint rule over each infinitesimal interval $\Delta t = 1/L$. There are many cosine expansions of the sinc function can be found from equation \eqref{A.3} by taking integral with help of efficient integration methods \cite{Mathews1999}. For example, another cosine expansions of the sinc function can be found by using the trapezoidal rule
\small
\[\label{A.4}
\tag{A.4}
{\rm{sinc}}\left( x \right) = \mathop {\lim }\limits_{L \to \infty } \frac{1}{L}\left[ {\frac{{1 + \cos \left( x \right)}}{2} + \sum\limits_{\ell  = 1}^{L - 1} {\cos \left( {\frac{\ell }{L}x} \right)} } \right]
\]
\normalsize
and the Simpson\text{'}s rule
\footnotesize
\[\label{A.5}
\tag{A.5}
{\rm{sinc}}\left( x \right) = \mathop {\lim }\limits_{L \to \infty } \frac{1}{{6L}}\left[ {1 + \cos \left( x \right) + 4\sum\limits_{\ell  = 1}^L {\cos \left( {\frac{{\ell  - 1/2}}{L}x} \right)}  + 2\sum\limits_{\ell  = 1}^{L - 1} {\cos \left( {\frac{\ell }{L}x} \right)} } \right].
\]
\normalsize
It is interesting to note that the limit \eqref{A.5} can also be derived trivially as a weighted sum of equations \eqref{A.2} and \eqref{A.4} in a proportion $2/3$ to $1/3$ as follows
\small
\[
\begin{aligned}
{\rm{sinc}}\left( x \right) = &\frac{2}{3} \times \mathop {\lim }\limits_{L \to \infty } \frac{1}{L}\sum\limits_{\ell  = 1}^{L} {\cos \left( {\frac{{\ell  - 1/2}}{L}x} \right)} \\
& + \frac{1}{3} \times \mathop {\lim }\limits_{L \to \infty } \frac{1}{L}\left[ {\frac{{1 + \cos \left( x \right)}}{2} + \sum\limits_{\ell  = 1}^{L - 1} {\cos \left( {\frac{\ell }{L}x} \right)} } \right].
\end{aligned}
\]
\normalsize

Any of these or similar cosine expansions of the sinc function can be used in integration to obtain expansion series for the error function ${\rm{erf}} \left( x \right)$ and, consequently, for the constant pi as well. However, as a simplest case we consider an application of equation \eqref{A.2} only. Thus, substituting the cosine expansion \eqref{A.2} of the sinc function into the integral \eqref{A.1} yields
\[
{\rm{erf}}\left( x \right) = \frac{{2x}}{\pi } \times \mathop {\lim }\limits_{L \to \infty } \int\limits_0^\infty  {\exp \left( { - {t^2}/4} \right)\underbrace {\frac{1}{L}\sum\limits_{\ell  = 1}^L {\cos \left( {\frac{{\ell  - 1/2}}{L}xt} \right)} }_{{\rm{sinc}}\left( {xt} \right)}dt}.
\]
Each terms in this equation is analytically integrable. Therefore, its integration leads to the expansion series \eqref{eq_3} of the error function. The more detailed description of the expansion series \eqref{eq_3} of the error function is given in our work \cite{Abrarov2016}.



\end{document}